\documentstyle[12pt]{article}
\parskip=\medskipamount                 
\thicklines                     
\newcommand{\bimn}[7]{\bibitem{#1}#2,
{\em #3},
{ #4}$\;${\bf
#5}$\;$(#6)$\;${#7}.}
\def\inbar{\vrule height1.5ex width.4pt depth0pt}
\def\IN{\relax{\rm I\kern-.18em N}}
\def\IQ{\relax\,\hbox{$\inbar\kern-.3em{\rm Q}$}}
\def\IR{\relax{\rm I\kern-.18em R}}
\def\ZZ{\relax{\sf Z\kern-.4em Z}}
\def\a{\alpha} \def\b{\beta}    
  
   \def\th{\theta}

 \def\cK{{\cal K}}  
 \def\cO{{\cal O}}  
 
\newtheorem{theorem}{Theorem}[section]

\newtheorem{corollary}{Corollary}[section]
\newtheorem{conjecture}{Conjecture}[section]
\newtheorem{lemma}{Lemma}[section]

\catcode`\@=11

\newif\if@fewtab\@fewtabtrue

\catcode`\@=11

\newif\if@fewtab\@fewtabtrue

{\count255=\time\divide\count255 by 60
\xdef\hourmin{\number\count255}
\multiply\count255 by-60\advance\count255 by\time
\xdef\hourmin{\hourmin:\ifnum\count255<10 0\fi\the\count255}}
\def\ps@draft{\let\@mkboth\@gobbletwo
    \def\@oddhead{}
    \def\@oddfoot
       {\hbox to 7 cm{\footnotesize {\em Draft version:} \draftdate

       \hfil}\hskip -7cm\hfil\rm\thepage \hfil}
    \def\@evenhead{}\let\@evenfoot\@oddfoot}

\def\ceqno{\global\@fewtabfalse
    \ifcase\@eqcnt \def\@tempa{& & &}\or \def\@tempa{& &}
      \or \def\@tempa{&}
      \or\def\@tempa{}\fi\@tempa
{\rm(\theequation)}}

\def\aeqno#1{\global\@fewtabfalse
    \ifcase\@eqcnt \def\@tempa{& & &}\or \def\@tempa{& &}
      \or \def\@tempa{&}
      \or\def\@tempa{}\fi\@tempa
{\rm(\theequation,#1)}}

\def\label#1{\ifnum\draftcontrol=1
 \global\def\draftnote{$\scriptstyle #1$}\fi
 \@bsphack\if@filesw {\let\thepage\relax
   \def\protect{\noexpand\noexpand\noexpand}%
\xdef\@gtempa{\write\@auxout{\string
      \newlabel{#1}{{\@currentlabel}{\thepage}}}}}\@gtempa
   \if@nobreak \ifvmode\nobreak\fi\fi\fi
  \@esphack}

\def\alabel#1#2{\label{#1}\global\@fewtabfalse
    \ifcase\@eqcnt \def\@tempa{& & &}\or \def\@tempa{& &}
      \or \def\@tempa{&}
      \or\def\@tempa{}\fi\@tempa
{\hbox to 3cm{\phantom{\rm(\theequation,#2)}
\draftnote \hfil}\hskip -3cm {\rm(\theequation,#2)}}}

\def\clabel#1{\label{#1}\global\@fewtabfalse
    \ifcase\@eqcnt \def\@tempa{& & &}\or \def\@tempa{& &}
      \or \def\@tempa{&}
      \or\def\@tempa{}\fi\@tempa
{\hbox to 3cm{\phantom{\rm(\theequation)}
\draftnote \hfil}\hskip -3cm{\rm(\theequation)}}}

\def\eqnarray{\def\draftnote{{}}\global\@fewtabtrue
\stepcounter{equation}\let\@currentlabel=\theequation
\global\@eqnswtrue
\global\@eqcnt\z@\tabskip\@centering\let\\=\@eqncr
$$\halign to \displaywidth\bgroup\@eqnsel\hskip\@centering\@eqcnt\z@
  $\displaystyle\tabskip\z@{##}$&\global\@eqcnt\@ne
  \hskip 1\arraycolsep \hfil$\displaystyle{##}$\hfil
  &\global\@eqcnt\tw@ \hskip 1\arraycolsep
$\displaystyle\tabskip\z@{##}$
\hfil  \tabskip\@centering&\global\@eqcnt\thr@@\llap{##}\tabskip\z@
\cr}

\def\endeqnarray{\@@eqncr\egroup
      \global\advance\c@equation\m@ne$$\global\@ignoretrue}

\def\@eqnnum{\hbox to 3cm{\phantom{\rm(\theequation)} \draftnote
                         \hfil}\hskip -3cm {\rm(\theequation)}}

\def\@@eqncr{\let\@tempa\relax
    \ifcase\@eqcnt \def\@tempa{& & &}\or \def\@tempa{& &}
      \or \def\@tempa{&}
      \or\def\@tempa{}
\fi\@tempa
\if@eqnsw
\if@fewtab\@eqnnum\fi
\stepcounter{equation}\fi\global
\@eqnswtrue\global\@eqcnt\z@\global\@fewtabtrue\cr}

\def\draftcite#1{\ifnum\draftcontrol=1#1\else{}\fi}

\def\@lbibitem[#1]#2{\item{}\hskip -3cm \hbox to 2cm
{\hfil$\scriptstyle\draftcite{#2}$}\hskip
1cm[\@biblabel{#1}]\if@filesw
     {\def\protect##1{\string ##1\space}\immediate
      \write\@auxout{\string\bibcite{#2}{#1}}}\fi\ignorespaces}

\def\@bibitem#1{\item\hskip -3cm \hbox to 2cm
{\hfil $\scriptstyle\draftcite{#1}$}\hskip 1cm
\if@filesw \immediate\write\@auxout
       {\string\bibcite{#1}{\the\value{\@listctr}}}\fi\ignorespaces}

\def\nsection#1{\section{#1}\setcounter{equation}{0}}

\def\draftdate{\number\month/\number\day/\number\year\ \ \ \hourmin }

\global\def\draftcontrol{0}
\catcode`\@=12

\def\theequation{{\thesection.\arabic{equation}}}

\def\qq{\begin{eqnarray}}
\def\qqq{\end{eqnarray}}
\def\rx#1{~(\ref{#1})}
\def\ex#1{eq.\rx{#1}}
\def\eex#1{eqs.\rx{#1}}
\def\cx#1{~\cite{#1}}

\newlength{\shiftwidth}
\def\shift#1{&&\hbox to \shiftwidth{\hfill $\displaystyle#1$}}
\newlength{\sshiftwidth}
\def\sshift#1{\lefteqn{\hbox to
\sshiftwidth{\hfill$\displaystyle#1$}}}

\def\ie{{\it i.e.\ }}
\def\eg{{\it e.g.\ }}

\def\rhs{{\it r.h.s.\ }}
\def\lhs{{\it l.h.s.\ }}

\def\p{^{\prime}}
\def\pp{^{\prime\prime}}

\def\prosign{\mathop{{\rm sign}}\nolimits}
\def\sign#1{\prosign\left(#1\right)}

\def\spint#1{\int\limits^{+\infty}_{\scriptstyle -\infty \atop
[{#1}]}}

\def\pr#1#2{ {\em Proof of #1~\ref{#2}.} }

\def\wc{ weak Conjecture~\ref{c1.1} }
\def\sc{ strong Conjecture~\ref{c1.2} }
\def\fone{ Fig.1 }

\def\jakk{ J_\a(\cK; K) }
\def\japq{ J_a(\kpq; K) }
\def\junk{ J_\a(\unknot; K) }
\def\vakk{ V_\a(\cK; K) }
\def\vapq{ V_\a(\kpq; K) }
\def\akz{ \nabla_A(\cK; z) }
\def\ankz{ \nabla_A^{2n+1}(\cK; z) }
\def\anpq{ \nabla_A^{2n+1}(\kpq; z) }
\def\apq#1{ \nabla_A\left( \kpq; #1 \right) }

\def\unknot{ {\rm unknot} }
\def\spka{ \sin \left( {\pi \over K} \a \right) }
\def\spk{ \sin \left( {\pi \over K} \right) }
\def\spkv#1{ \sin \left( {\pi \over K} #1 \right) }

\def\fu#1{ \Phi \left(U^{(#1)} \right) }
\def\ut#1#2{ \tilde{U}^{(#1)}_{#2} }

\def\ipf{ {i\pi\over 4} }
\def\sqtk{ \sqrt{2\over K} }
\def\sqkt{ \sqrt{K\over 2} }
\def\iptk{ {i\pi\over 2K} }
\def\pqf{ \left( \xp\xq-{\xp\over\xq}-{\xq\over\xp} \right) }

\def\spq{ \check{q} }
\def\th{ \tilde{\xh} }
\def\xp{ p }
\def\xq{ q }
\def\xr{ r }
\def\xs{ s }
\def\xt{ t }
\def\xh{ h }

\def\qat{ \spq^{\a \over 2} }
\def\qatp{ \left( \spq^{\a \over 2} \right) }
\def\qatm{ \spq^{- {\a \over 2} } }
\def\qt{ \spq^{1\over 2} }
\def\qtm{ \spq^{-{1\over 2} } }
\def\qad{ \qat - \qatm }

\def\tpik{ {2\pi i \over K} }
\def\iptk{ {i\pi \over 2K} }
\def\ipk{ {i\pi \over K} }

\def\smnz{ \sum_{m,n\geq 0} }

\def\vnkz{ V^{(n)} (\cK;z) }
\def\vzkz{ V^{(0)} (\cK;z) }
\def\vnpq{ V^{(n)} (\kpq;z) }
\def\tvkz#1{ \tilde{V}^{(#1)}(\cK;z) }
\def\pnkz{ P^{(n)} (\cK;z) }
\def\pnpq{ P^{(n)} (\kpq;z) }
\def\tpz#1{ \tilde{P}^{(#1)}(\cK;z) }
\def\pappr#1{ P^{(#1)}_{\rm appr} }
\def\tpappr#1{ \tilde{P}^{(#1)}_{\rm appr} }

\def\DmnK{ D_{m,n}(\cK) }
\def\dnm{ d^{(n)}_m }
\def\tdm#1{ \tilde{d}^{(#1)}_m }

\def\kpq{ \cK_{\xp,\xq} }

\def\tm#1{ \left( \xt^{#1} - \xt^{-#1} \right) }
\def\tmo{ \xt - \xt^{-1} }

\def\sk#1{ \sum_{1\leq #1 \leq K-1} }

\def\mexp#1{ \left( e^{i\pi#1 {\a_\b \over K} }
- e^{ - i\pi#1 {\a_\b \over K} } \right) }

\def\zz#1{ \ZZ\left[ #1 \right] }
\def\zzs#1{ \ZZ[ #1 ] }
\def\zzf#1{ \ZZ\left[ z^2, {1\over #1} \right] }

\def\mpq{ \mu\xp - \xq - \xp\xq }
\def\pqm{ {1\over 2}(\xp + 1)(\xq - \mu) }
\def\qcas{ \spq^{ {1\over 2} (\xp+1)(\xq -\mu) } }

\def\qg{ \spq^{-{1\over 4} {\b^{\prime 2}\over \xp\xq} } }
\def\qfr{ {1-\spq^{\xq\a + {\b\over \xp} }\over 1-\spq^{\xp\xq\a+\b}
} }
\def\zset{ (\ZZ_+) - \left( \ZZ_+ +{1\over \xp} \right) }

\def\cmkpq{ C_{m,k}(\xp,\xq) }
\def\thee{ T(\xh, \epsilon) }
\def\tphe{ T^p(\xh, \epsilon) }

\def\zqa{ \ZZ \left[ \qat, \qatm \right] }
\def\zqap{ \ZZ \left[ \qat, \qatm, {1\over \xp} \right] }

\def\qpm{ \left( \qat, \qatm \right) }
\def\qfrp{ \left( {1-\qatp^{2\xp\xq} \over 1-\qatp^{2\xq} } \right)}

\def\vmin{ v_{\rm min} }

\begin{document}

\begin{titlepage}
\centerline{\hfill                 q-alg/9601009}
\vfill
\begin{center}

{\large \bf
Higher Order Terms in the Melvin-Morton Expansion of the Colored
Jones Polynomial.
}
\\

\bigskip
\centerline{L. Rozansky
}

\centerline{\em School of Mathematics, Institute for Advanced Study}
\centerline{\em Princeton, NJ 08540, U.S.A.}
\centerline{{\em E-mail address: rozansky@math.ias.edu}}

\vfill
{\bf Abstract}

\end{center}
\begin{quotation}
We formulate a conjecture about the structure of `upper lines'
in the expansion of the colored Jones polynomial of a knot in
powers of $(q-1)$. The Melvin-Morton conjecture states that the
bottom line in this expansion is equal to the inverse Alexander
polynomial of the knot. We conjecture that the upper lines are
rational functions whose denominators are powers of the Alexander
polynomial. We prove this conjecture for torus knots and give
experimental evidence that it is also true for other types of knots.
\end{quotation}
\vfill
\end{titlepage}

\pagebreak

\nsection{Introduction}
\label{s1}
Ever since the discovery of the Jones polynomial, its relation to the
objects of the classical topology, \ie the fundamental group of a
knot, remained somewhat of a mystery. An apparent similarity between
the skein relations for the Jones and Alexander polynomials did not
lead to a better understanding of this relation. Therefore the
discovery by P.~Melvin and H.~Morton\cx{MM} of the inverse Alexander
polynomial inside the $(q-1)$ expansion of the colored Jones
polynomial was a very interesting development.

Let $\cK$ be a knot in $S^3$. We denote by $\jakk$ its colored Jones
polynomial normalized in such a way that it is multiplicative under a
disconnected sum and
\qq
\junk = {\spka \over \spk} = {\qat - \qatm \over \qt - \qtm}, \qquad
\spq = e^\tpik.
\label{1.1}
\qqq
Another polular normalization for the Jones polynomial is
\qq
\vakk = {\jakk \over \junk}, \qquad
V_\a \in \ZZ[\spq, \spq^{-1}].
\label{1.2}
\qqq

For a fixed value of color $\a$ we can expand the Jones polynomial
$\vakk$ in Taylor series in powers of
\qq
\xh = \spq - 1,
\label{1.3}
\qqq
or, equivalently, in powers of
\qq
{1\over K} = {1\over 2\pi i} \log (1 + \xh).
\label{1.03}
\qqq
The coefficients of this expansion are polynomials of finite degree
in $\a$:
\qq
\vakk = \smnz \DmnK \a^{2m} \xh^n, \qquad
\DmnK \in \IQ
\label{1.4}
\qqq
The coefficients $\DmnK$ are rational invariants of the knot $\cK$.
D.~Bar-Natan\cx{BN} and J.~Birman, X-S.~Lin\cx{BL} showed that
$\DmnK$ were Vassiliev invariants of order $n$.

The following theorem was conjectured by P.~Melvin and
H.~Morton\cx{MM} and later proved by D.~Bar-Natan and
S.~Garoufalidis\cx{BG} (for a simple path integral proof
see\cx{Ro1}).
\begin{theorem}
\label{t1.1}
Let $\cK$ be a knot in $S^3$. Then the coefficients $\DmnK$ of the
expansion\rx{1.4} satisfy the following two properties
\qq
&
\DmnK = 0 \qquad for \qquad m\geq {n\over 2},
\label{1.5}
\\
&
\sum_{m\geq 0} D_{m, 2m}(\cK) a^{2m} =
{1 \over \nabla_A(\cK; e^{i\pi a} - e^{- i\pi a}) },
\label{1.6}
\qqq
here $\akz$ is the Alexander-Conway polynomial satisfying the skein
relation ``$(X) - (X) = z(II)$'' and normalized is such a way that
\qq
\nabla_A(\unknot;z) = 1.
\label{1.7}
\qqq
\end{theorem}

The bound on the powers of $\a$ in the expansion\rx{1.4} allows us to
rearrange it in `lines'
\qq
\vakk = \sum_{n\geq 0} \xh^n \sum_{m\geq 0} D_{m, n+2m} (\a\xh)^{2m}.
\label{1.8}
\qqq
From the quantum field theory point of view, the $n$th line
$\sum_{m\geq 0} D_{m,(n-1) + 2m} (\a\xh)^{2m}$ is related to the
$n$-loop contribution in the calculation of $\vakk$ as a Chern-Simons
path integral over the $SU(2)$ connections in the knot complement
(see\cx{Ro1} for details).

The Taylor expansions\rx{1.4}, \rx{1.8} and their link analogs played
a key role in defining the `perturbative' invariants of rational
homology spheres
(see a review in\cx{Ro6} and references therein)
and in establishing their relation\cx{Ro5} to Ohtsuki's
invariants\cx{Oh2}. The latter application prompted us to look for
`integrality' properties of the coefficients $\DmnK$. We conjecture
that this integrality can be exposed by an appropriate choice of
expansion parameters in the series\rx{1.4} and \rx{1.8}.

Let us introduce a new variable
\qq
z=\qad.
\label{1.9}
\qqq
%
%
%
We can express $\a\xh$ as
a power series in $z$ and $\xh$ by expanding
the \rhs of the equation
\qq
\a\xh = 2\log\left(\sqrt{1 + \left( {z\over 2} \right)^2 } + {z\over
2} \right) {h\over \log(1+h)} = z + \cO(z^3,\xh).
\label{1.10}
\qqq
After putting this expression in place of $(\a\xh)$ in \ex{1.8} and
assembling the powers of $\xh$ and $z$ we get a new expansion of the
Jones polynomial
\qq
\vakk & = & \sum_{n=0}^\infty \vnkz \xh^n,
\label{1.11}
\\
\vnkz & = & \sum_{m=0}^\infty \dnm(\cK) z^{2m}.
\label{1.12}
\qqq
The form of the substitution\rx{1.10} suggests immediately that the
bottom line in the expansion\rx{1.8} does not change:
\qq
\sum_{m=0}^\infty D_{m,n+2m} (\a\xh)^{2m} = \vzkz + \cO(\xh).
\label{1.13}
\qqq
As a result, the second part of the Melvin-Morton conjecture\rx{1.6}
takes the form
\qq
\vzkz = {1\over \akz}
\label{1.14}
\qqq
in new variables $z,\xh$. Since $\akz \in \zzs{z^2}$ and
$\nabla_A(\cK;0) = 1$, it follows from \ex{1.14} that
\qq
d_m^{(0)} \in \ZZ.
\label{1.15}
\qqq
We conjecture that the upper lines $\vnkz$ satisfy the properties
similar to those of\rx{1.14} and\rx{1.15}.
\begin{conjecture}[weak]
\label{c1.1}
All coefficients $\dnm$ in the expansion of the colored Jones
polynomial $\vakk$ of a knot $\cK\in S^3$ in powers  of $z = \qad$
and $\xh = \spq - 1$ are integer:
\qq
\dnm \in \ZZ, \qquad m,n\geq 0.
\label{1.16}
\qqq
\end{conjecture}
\begin{conjecture}[strong]
\label{c1.2}
A line $\vnkz$ in the series\rx{1.11} is a rational function of $z$:
\qq
\vnkz = {\pnkz \over \ankz}, \qquad \pnkz \in \zzs{z^2}.
\label{1.17}
\qqq
In other words, the expansion of the \rhs of \ex{1.17} in powers of
$z$ produces the series in the \rhs of \ex{1.12}.
\end{conjecture}
The weak conjecture can be derived from the strong one in the same
way in which we derived\rx{1.15} from \ex{1.14}.

It may happen that for some knots the polynomials $\pnkz$ defined by
\ex{1.17} would be divisible by powers of $\akz$. This means that for
those knots a smaller power of $\akz$ could be placed in
the denominator.
Amphicheiral knots seem to present an example of such behavior.
These are the knots which are isotopic to
their mirror image.

If $\cK\p$ is the mirror image of $\cK$, then
\qq
V_\a(\cK\p;K) = V_\a(\cK;-K).
\label{3.2}
\qqq
This symmetry is not easily seen in the coefficients of
expansion\rx{1.11} because it transforms $\xh$ into
$-{\xh\over 1+\xh}$ rather than into $-\xh$. Therefore it is natural
to try another expansion parameter
\qq
\th =  \qt - \qtm = (1+\xh)^{1\over 2} - (1+\xh)^{-{1\over 2}}
\label{3.3}
\qqq
instead of $\xh$ in \ex{1.11}:
\qq
\vakk & = &
\sum_{n=0}^{\infty} \tvkz{n} \th^n,
\label{3.4}
\\
\tvkz{n} & = & \sum_{m=0}^{\infty} \tdm{n}(\cK) z^{2m}.
\label{3.5}
\qqq
The symmetry\rx{3.2} converts $\th$ into $-\th$
and $z$ into $-z$,
so for amphicheiral
knots
\qq
\tvkz{2n+1} = 0, \qquad n\geq 0.
\label{3.6}
\qqq

Since the relation between $\xh$ and $\th$ involves fractional
powers, it does not follow from the \wc that the coefficients
$\tdm{n}$ would also be integer
for any knot. In fact, our numerical estimates
show that some of the first coefficients for the knot $6_1$ are
fractional. However, for the amphicheiral knots only the even
powers of $\th$ participate in the expansion\rx{3.4} due to \ex{3.6}.
Since the expansion
\qq
\th^2 = \sum_{n=2}^{\infty} (-1)^n \xh^n
\label{3.06}
\qqq
contains only integer coefficients (\ie $(-1)^n$), then the \wc
implies that for amphicheiral knots, $\tdm{n}$ should also be integer:
\begin{corollary}
\label{cor3.1}
For an amphicheiral knot $\cK$ in addition to \rx{1.16}
\qq
\tdm{n} \in \ZZ.
\label{3.7}
\qqq
\end{corollary}
Our experimental data also suggests (see Section~\ref{s3} for
details) that the following enhancement of the \sc is true
\begin{conjecture}
\label{c3.2}
For an amphicheiral knot $\cK$ a line $\tvkz{2n}$ in the
series\rx{3.4} is a rational function
\qq
\tvkz{2n} = {\tpz{2n} \over \nabla_A^{3n+1}(\cK;z) }, \qquad
\tpz{2n} \in \zzs{z^2}.
\label{3.8}
\qqq
\end{conjecture}
Note that $3n+1\leq 2(2n) + 1$, $2(2n)+1$ being the power required by
\ex{1.17}, so the amphicheiral knots require a smaller power of the
Alexander polynomial in denominators.

The integrality of the coefficients of the polynomials $\pnkz$ gives
us a hope that similarly to the denominator
$\ankz$, they may also have a direct
interpretation in the framework of classical topology.

The \wc has a `practical' application: we will
use it in\cx{Ro8} in order to derive a $p$-adic convergence of the
series of perturbative invariants to the total
Witten-Reshetikhin-Turaev invariant of rational homology spheres
constructed by rational surgeries on a knot in $S^3$.

In Section~\ref{s2} we derive the \sc for torus knots. In
Section~\ref{s3} we present experimental evidence that our
conjectures are also true for other types of knots. In
Section~\ref{s4} we speculate about the possible explanation for the
power of the Alexander polynomial in denominators of \eex{1.17}
and\ex{3.8}.

\nsection{The Jones polynomial of torus knots}
\label{s2}
We denote a type $(\xp,\xq)$ torus knot in $S^3$ as $\kpq$. An
expression for its colored Jones polynomial was derived in\cx{ILR}
within the framework of the quantum Chern-Simons theory. The $(q-1)$
expansion of the polynomial was studied in\cx{M} and\cx{AL}. We
derived the explicit expansion\rx{1.8} for torus knots in
eq.(A.4) of Appendix of\cx{Ro1}. In our notations there
\qq
Z_\a(S^3, \cK_{m,n};K-2) = \sqrt{2\over K} \spka V_\a(\cK_{m,n};K).
\label{2.1}
\qqq
We will rederive the formula
of\cx{Ro1} in a slightly different way that will
allow us to prove the conjectures of the previous section.
\begin{lemma}
\label{l2.1}
The expansion\rx{1.11} for a torus knot $\kpq$ comes from the formula
\qq
\vapq & = &
{1\over z} { e^{ \ipf\sign{\xp\xq} } \over \sqrt{2K|\xp\xq|} }
e^{\iptk \pqf}
\spint{\b=0} d\b\,e^{-\iptk {\b^2\over \xp\xq} }
{z_\b \over \apq{z_\b} }
\label{2.17}
\\
& = &
\!\!\!
\left.
{1\over z}
(1+\xh)^{ {1\over 4} \pqf }
\!\!\!
\sum_{m=0}^{\infty} {(2m)!\over (m!)^2}
\!\!
\left( {\log(1+\xh) \over 4\xp\xq} \right)^m
\!\!\!
\left( {K\xp\xq \over i\pi} \right)^{2m}
{\partial^{(2m)} \over \partial\b^{2m} }
{z_\b \over \apq{z_\b} } \right|_{\b=0},
\nonumber
\qqq
here
\qq
z_\b = \spq^{ {1\over 2} \left(\a + {\b\over \xp\xq} \right) } -
\spq^{ -{1\over 2} \left(\a + {\b\over \xp\xq} \right) } ,
\label{2.1017}
\qqq
$\apq{z}$ is the Alexander polynomial of $\kpq$:
\qq
\apq{\tmo} = {\tm{\xp\xq}\left(\tmo\right) \over \tm{\xp} \tm{\xq} }.
\label{2.15}
\qqq
and the symbol $\spint{\b=0}$ means that we have to take only the
contribution of the stationary phase point $\b=0$ to the integral of
\ex{2.12}.
\end{lemma}
Note that each derivative $\partial_\b$ extracts a factor of
$i\pi\over K\xp\xq$ from $z_\b$. These factors cancel the
prefactor $\left( {K\xp\xq \over i\pi} \right)^{2m}$, so that
\ex{2.17} presents an expansion in powers of
\qq
\log(1+\xh) = \xh - \sum_{n=2}^{\infty}(-1)^n {\xh^n \over n}.
\label{2.017}
\qqq
\pr{Lemma}{l2.1}
If $\kpq$ is a torus knot, then the numbers $\xp,\xq\in\ZZ$ are
coprime. Therefore we can choose the numbers $\xr,\xs\in\ZZ$ in such
a way that
\qq
\xp\xs -\xq\xr = 1.
\label{2.2}
\qqq
It is not hard to see that the surgeries on the 3-component solid
link in $S^3$ of \fone (with framings
$\left(-{p\over r}, {q\over s}, 0 \right)$) produce again $S^3$.
However the dashed unknot of \fone becomes a torus knot $\kpq$ in the
new $S^3$. Therefore its colored Jones polynomial can be calculated
by a Reshetikhin-Turaev surgery formula applied to the link of \fone.
The colored Jones polynomial of that link is
\qq
J_{\a_1,\a_2,\a,\b} =
\frac{ \spkv{\a_1\b} \spkv{\a_2\b} \spkv{\a\b} }
{\sin^2\left( {\pi\over K} \b \right) \spk }.
\label{2.3}
\qqq
The `quantum factor' for a rational surgery with the framing
$p\over q$ was worked out by L.~Jeffrey\cx{Je}:
\qq
\ut{\xp,\xq}{\a 1}
& = &
i {\sign{q}\over \sqrt{2K|q|} }e^{-\ipf\fu{\xp,\xq} }
\sum_{n=0}^{q-1} \sum_{\mu=\pm 1} \mu
\label{2.4}
\\
&&
\qquad\times
\exp\left( \iptk {1\over \xq}
(p\a^2 - 2\a(2Kn+\mu) + s(2Kn+\mu)^2 )
\right),
\nonumber
\qqq
here $\fu{\xp,\xq}$ is the Rademacher function:
\qq
\Phi\left[
\matrix{\xp & \xr \cr \xq & \xs\cr}
\right] =
{\xp + \xs \over \xq} - 12 s(\xp,\xq),
\label{2.5}
\qqq
and $s(\xp,\xq)$ is the Dedekind sum:
\qq
s(\xp,\xq) = {1\over 4\xq}
\sum_{j=1}^{q-1} \cot\left(\pi {j\over \xq} \right)
\cot \left( \pi {\xp j\over \xq} \right).
\label{2.6}
\qqq
Since
\qq
\sk{\a} \spkv{\a\b}\ut{\xp,\xq}{\a 1} =
\sqkt\ut{-\xq,\xp}{\b 1},
\label{2.7}
\qqq
we conclude that
\qq
\japq = { e^{-\iptk \xp\xq(\a^2 - 1) } \over \sqtk \spk}
\sk{\b} {\spkv{\a\b}\over\spkv{\b} }
\ut{\xr,\xp}{\b 1}\ut{-\xs,\xq}{\b1}.
\label{2.8}
\qqq
We did not include the manifold framing correction because the
surgery produces $S^3$ in the canonical framing. However we had to
include the knot framing correction factor
$e^{-\iptk \xp\xq(\a^2 - 1) }$ because the torus knot is produced
with the framing $\xp,\xq$.

By substituting \ex{2.4} in \ex{2.8} and using the relation
\qq
\Phi \left[ \matrix{ \xr & -\xs \cr \xp & -\xq \cr} \right] +
\Phi \left[ \matrix{ -\xs & \xr \cr \xq & -\xp \cr} \right] =
- 3\sign{\xp\xq},
\label{2.9}
\qqq
we arrive at the following formula for the colored Jones polynomial:
\qq
\japq & = &
\frac { e^{ \ipf\sign{\xp\xq} } } { 4\sqrt{2K|\xp\xq|} }
\frac { e^{-\iptk \xp\xq (\a^2-1)} } {\spk}
\sum_{n_1=0}^{p-1} \sum_{n_2=0}^{q-1} \sum_{\mu_{1,2,3}=\pm 1}
\mu_1\mu_2\mu_3
\nonumber
\\
&&\times
\!\!\!\!\sk{\b}
{1\over\spkv{\b}} \exp \left[ -\iptk \left( {\b^2\over \xp\xq}
+ 2\b \left( {2Kn_1+\mu_1 \over \xp} + {2Kn_2+\mu_2 \over \xq} +
\mu_3\a \right)
\right.\right.
\nonumber\\
&&\qquad \left.\left.
+ {\xq \over \xp} (2Kn_1+\mu_1)^2 + {\xp \over \xq} (2Kn_2+\mu_2)^2
\right) \right].
\label{2.10}
\qqq
The sum over $\b$ in this formula is completely similar to the sums
in eqs.~(2.8), (2.9) of \cx{RoS2} if we substitute there
$g=0,\; n=3,\; m_1=n_1,\; m_2=n_2,\;K{m_3\over p_3}=\mu_3\a$. We will
not present the analysis of the large $K$ asymptotics of the
formula\rx{2.10} since it is exactly the same as the one in Section~3
of\cx{RoS2}. We will rather use the final result expressed in
eqs.~(3.21), (3.23) and eq.~(3.44) of the Proposition~3.1 in that
paper. Namely, if
\qq
\a < {K\over |\xp\xq|},
\label{2.11}
\qqq
then all the `irreducible' contributions~(3.45) of\cx{RoS2} cancel
out (note, that there are no irreducible flat connections in the knot
complement whose holonomy along the knot meridian is equal to
$\exp\left(\ipk \a\sigma_3\right)$,
$\sigma_3 = \left( \matrix{ 1 & 0 \cr 0 & -1 \cr } \right) $ for $\a$
satisfying\rx{2.11} ). The only survivor is the `reducible'
contribution
\qq
\japq & = &
\frac { e^{ \ipf\sign{\xp\xq} } } { 8\sqrt{2K|\xp\xq|} }
\frac { e^{-\iptk \xp\xq (\a^2-1)} } {\spk}
\sum_{\mu_{1,2,3}=\pm 1}
\mu_1\mu_2\mu_3
\label{2.12}
\\
&&\times
\spint{\b = - \xp\xq\mu_3\a}
{d\b \over \spkv{\b} } \exp \left[ -\iptk \left( {\b^2\over \xp\xq} +
2\b \left( {\mu_1\over \xp} + {\mu_2\over \xq} + \mu_3\a \right) +
{\xq\over \xp} + {\xp \over \xq} \right) \right]
\nonumber\\
& = &
\frac { e^{ \ipf\sign{\xp\xq} } } { \sqrt{2K|\xp\xq|} }
\frac { e^{ \iptk \pqf } }
{2i\spk}
\nonumber
\\
&&\times
\spint{\b=0} d\b\, e^{-\iptk {\b^2\over \xp\xq} }
{\mexp{\xp}\mexp{\xq} \over \mexp{\xp\xq} },
\nonumber
\qqq
here we used a notation
\qq
\a_\b = \a + {\b \over \xp\xq}.
\label{2.13}
\qqq
The symbol $\spint{\b=0}$ means that we have to take only the
contribution of the stationary phase point $\b=0$ to the integral of
\ex{2.12}. In other words, we have to expand the preexponential
factor in Taylor series in $\b$ around $\b=0$ and integrate each
monomial with the gaussian factor
$e^{ -\iptk {\b^2\over \xp\xq} }$ term by term with the help of the
formula
\qq
\int_{-\infty}^{+\infty} e^{-\iptk {\b^2\over \xp\xq} } \b^{2m} d\b =
\sqrt{2K|\xp\xq|} e^{-\ipf\sign{\xp\xq} } {(2m)!\over m!}
\left( {K\xp\xq \over 2\pi i} \right)^m.
\label{2.14}
\qqq
The result will be precisely eq.~(A.4) of\cx{Ro1}. We present
that formula in a slightly different form. We use the
formula\rx{2.15} for the Alexander polynomial of the torus knot
and by introducing a notation
\qq
z_\b = \spq^{\a_\b\over 2} - \spq^{-\a_\b\over 2}
\label{2.16}
\qqq
we arrive at \ex{2.17}. $\Box$

The formula\rx{2.17} proves half of the \sc, namely,
\begin{lemma}
\label{l2.2}
For a torus knot $\kpq$,
\qq
\vnpq = {\pnpq \over \anpq}, \qquad
\pnpq \in \IQ[z^2]
\label{2.019}
\qqq
\end{lemma}
\pr{Lemma}{l2.1}
For a smooth function $f(z_\b)$
\qq
\left( {K\xp\xq \over i\pi} \right)^2
{\partial^2\over \partial \b^2} f(z_\b) =
z_\b f\p(z_\b) + (z_\b^2 + 4) f\pp(z_\b),
\label{2.18}
\qqq
therefore
\qq
\left.
\left( {K\xp\xq\over i\pi} \right)^{2m}
{\partial^{(2m)} \over (\partial\b)^{2m} }
{z_\b \over \apq{z_\b} } \right|_{\b=0}
=
{z \tilde{P}_m(z) \over \nabla_A^{2m+1}(\kpq;z) },
\qquad \tilde{P}_m(z) \in \zzs{z^2}.
\label{2.19}
\qqq
The numerator of the \rhs of \ex{2.19} is proportional to $z$ because
in view of \ex{2.18}, the \lhs of this equation is an odd function
of $z$. Finally, \ex{2.17} demonstrates that a line $\vnpq$ is a
linear combination of the functions\rx{2.19} for $m\leq n$. $\Box$


The fractions in the expansion of $(1+\xh)^{ {1\over 4}\pqf }$ and
$\log(1+\xh)\over 4\xp\xq$ in powers of $\xh$ separate us from the
complete proof of the \sc. We prove the next lemma in order to
eliminate half of these fractions.
\begin{lemma}
\label{l2.3}
For a torus knot $\kpq$,
\qq
\pnpq \in \zzf{\xp}
\label{2.1019}
\qqq
\end{lemma}
\pr{Lemma}{l2.3}
We are going to absorb the
factor
$$\mexp{\xp}$$
together with some other factors of the integrand in \ex{2.12}
inside the
gaussian factor
$e^{-\iptk {\b^2\over \xp\xq} }$ by completing the square. We achieve
this by introducing a new variable
$\b\p$ such that
\qq
\b = \b\p+\mpq,
\label{2.20}
\qqq
Now the integral of \ex{2.12} can be rewritten as
\qq
\vapq  =
{1\over z} {e^{\ipf\sign{\xp\xq} }\over \sqrt{2K|\xp\xq|} }
\sum_{\mu=\pm 1} \qcas \spq^{ {\a\over 2}(\mpq)}
\spint{\b\p=0} d\b\p\,\qg \qfr.
\label{2.21}
\qqq
We kept $\b$ in the last factor of the integrand in this equation
meaning that it is a function\rx{2.20} of $\b\p$.

The last factor of the integrand in \ex{2.21} can be
presented as a geometric series
\qq
\qfr = \lim_{x\rightarrow 1^-} \sum_{n\in\zset}
\left( \spq^{\xp\xq\a+\b}\right)^n x^n,
\label{2.22}
\qqq
here we used a notation
\qq
\sum_{n\in\zset} f(n) = \sum_{n=0}^{\infty} f(n) -
\sum_{n=0}^{\infty} f(n+{1\over p}).
\label{2.23}
\qqq
Since
\qq
{e^{\ipf \sign{\xp\xq} }\over \sqrt{2K|\xp\xq|} }
\int_{-\infty}^{+\infty} d\b\p\,\qg \spq^{n\b\p} =
\spq^{-\xp\xq n^2},
\label{2.24}
\qqq
\ex{2.21} becomes
\qq
\vapq & = & {1\over z} \sum_{\mu=\pm 1} \mu (1+\xh)^{\pqm}
\left(\qat\right)^{\mpq}
\label{2.25}
\\
&&
\qquad\times
\lim_{x\rightarrow 1^-} \sum_{n\in\zset}
\left(\qat\right)^{2\xp\xq n} (1+\xh)^{n(\mpq)} (1+\xh)^{-\xp\xq n^2}
x^n.
\nonumber\\
\qqq
Our immediate task is to expand this expression in powers of $\xh$
with coefficients being hopefully rational functions of $\qat$. Since
$\xp$ and $\xq$ are coprime, at least one of them is odd. Therefore
$\pqm \in \ZZ$ and
\qq
(1+\xh)^{\pqm} \in \ZZ[[\xh]].
\label{2.26}
\qqq
It remains to treat the sum over $n$. Consider an expansion
\qq
(1+\xh)^{-\xp\xq n^2} = \sum_{m=0}^{\infty} {(-1)^m}
{\prod_{l=0}^{m-1} (\xp\xq n^2+l) \over m!} \xh^m .
\label{2.27}
\qqq
Each term in the sum over $m$ is a polynomial in $n$ which takes
integer values for $n\in \ZZ$. Therefore it can be presented as an
integer linear combination of binomial polynomials
${n \choose k} = {n! \over (n-k)! k!}$:
\qq
{(-1)^m} {\prod_{l=0}^{m-1} (\xp\xq n^2+l) \over m!}=
\sum_{k=0}^{2m} \cmkpq {n \choose k}, \qquad
\cmkpq \in \ZZ.
\label{2.28}
\qqq
A binomial polynomial can be expressed as a derivative
\qq
\left.
{n \choose k} = {1\over k!} \partial^{(k)}_\epsilon (1+\epsilon)^n
\right|_{\epsilon=0}
\label{2.29}
\qqq

We combine \eex{2.27} -\rx{2.29} together and substitute them back
into \ex{2.25}. After summing up a geometric series
\qq
&
\lim_{x\rightarrow 1^-} \sum_{n\in \zset} \qatp^{2\xp\xq n}
(1+\xh)^{n(\mpq)} (1+\epsilon)^n x^n =
{1 - \thee \over 1- \tphe},
&
\label{2.30}
\\
&
\thee = \qatp^{2\xq} (1+\xh)^{\mu-\xq-{\xq\over\xp} }
(1+\epsilon)^{1\over p},
\label{2.32}
&
\qqq
we get
\qq
\vapq & = & {1\over z} \sum_{\mu=\pm 1} \mu (1+\xh)^{\pqm}
\left(\qat\right)^{\mpq}
\label{2.31}
\\
&&
\qquad\times\left.
\sum_{m=0}^{\infty} \xh^m \sum_{k=0}^{2m} \cmkpq
\left( {1\over k!} \partial_\epsilon^{(k)}\right)
{1 - \thee \over 1- \tphe} \right|_{\epsilon=0}.
\nonumber
\qqq
We see from this expression that a line $\vnpq$ of the
expansion\rx{1.11} is a linear combination over $\zqa$ of
coefficients at monomials
\qq
\xh^{n_1} \epsilon^{n_2}, \qquad
n_1 + {n_2\over 2} \leq n
\qqq
in the expansion of
${1 - \thee \over 1- \tphe}$. A coefficient at
$\xh^{n_1} \epsilon^{n_2}$ would have taken the form
\qq
{ P_{n_1,n_2}\qpm \over \qfrp^{n_1+n_2+1} }
\label{2.34}
\qqq
with some polynomial
$P_{n_1,n_2} \in \zqa$ if not for the fractional powers in
$(1+\xh)^{1\over \xp}$ and $(1+\epsilon)^{1\over \xp}$. However a
simple lemma
\qq
{1\over n!} \prod_{l=0}^{n-1} \left({1\over \xp} + l \right)
\in \zz{{1\over \xp}}
\label{2.35}
\qqq
guarantees that
\qq
P_{n_1,n_2} \in \zqap.
\label{2.36}
\qqq
Therefore we find that
\qq
\vnpq = {P_n\qatp \over \qfrp^{2n+1} }, \qquad
P_n\qatp \in \zqap.
\label{2.37}
\qqq
By using \ex{2.15} we can rewrite this as
\qq
&
\vnpq =
\frac{\tilde{P}_n\qatp}
{ \left( 1 - \qatp^{2\xp} \right)^{2n+1} \nabla_A^{2n+1}
\left( \kpq; \qat-\qatm \right) },
&
\label{2.38}
\\
&
\tilde{P}_n\qatp \in \zqap.
&
\nonumber
\qqq
Comparing this with \ex{2.019} we conclude that
$\tilde{P}_n\qatp$ is divisible by
$\left(1 - \qatp^{2\xp}\right)^{2n+1}$
over $\IQ\left[ \qat, \qatm \right]$. However
since the process of division by $1 - \qatp^{2\xp}$ does not
introduce any new fractions in the coefficients of the polynomials,
the polynomial $\pnpq$ of \ex{2.019} can have only
the divisors of $\xp$ in denominators of its coefficients:
\qq
P^{(n)} \left( \kpq; \qat - \qatm \right) =
{\tilde{P}_n\qatp \over
\left( 1 - \qatp^{2\xp} \right)^{2n+1} } \in \zqap.
\label{2.39}
\qqq
Since, according to Lemma~\ref{l2.2}, $\pnpq \in \IQ[z^2]$, \ex{2.39}
proves the Lemma. $\Box$
%
%

\noindent
\pr{Conjecture}{c1.2}
A proof similar to that of Lemma~\ref{l2.3} would prove that
$$\pnpq \in \zzf{\xq}.$$
Since
$\xp$ and $\xq$ are coprime, $\zzf{\xp} \bigcap \zzf{\xq}=\zzs{z^2}$
and
this proves the \sc. $\Box$


We used the formula\rx{2.17} to calculate the polynomials
$P^{(n)}(z)$, $n=1,2,3$ for the simplest torus knots:
%
\qq
\begin{array}{lrcl}
(2,3): \qquad &
\nabla_A & = & 1 + z^2 \hspace{298.pt}
\\
& P^{(1)}& = & 2z^2 + z^4
\\
& P^{(2)}& = & 1 - 3z^2 - z^3
\\
& P^{(3)}& = & -3 + 13z^2 - z^6
\\
\end{array}
\label{e2.1}
\qqq
\qq
\begin{array}{lrcl}
(2,5): \qquad  &
\nabla_A & = & 1 + 3z^2 + z^4 \hspace{266.7pt}
\\
& P^{(1)}& = & 10z^2 + 21z^4 + 12z^6 + 2z^8
\\
& P^{(2)}& = & 3 - 19z^2 - 24z^4 + 58z^6 + 145z^8 + 128z^{10} +
56z^{12} + 12z^{14} + z^{16}
\end{array}
\label{e2.2}
\qqq
\qq
\begin{array}{lrcl}
(2,7): \qquad  &
\nabla_A & = & 1 + 6z^2 + 5z^4 + z^6  \hspace{235.6pt}
\\
& P^{(1)}& = & 28z^2 + 126z^4 + 180z^6 + 110z^8 + 30z^{10} + 3z^{12}
\\
& P^{(2)}& = & 6 - 66z^2 - 138z^4 + 1398z^6 + 7248z^8 + 15747z^{10} +
19635z^{12}
\\
&&&
\qquad
+ 15360z^{14} + 7776z^{16} + 2544z^{18} + 519z^{20} +
60z^{22} + 3z^{24}
\end{array}
\label{e2.3}
\qqq
\qq
\begin{array}{lrcl}
(3,5): \qquad &
\nabla_A& = & 1 + 8z^2 + 14z^4 + 7z^6 + z^8
\\
& P^{(1)}& = & 40z^2 + 314z^4 + 908z^6 + 1224z^8 + 846z^{10} +
308z^{12} + 56z^{14} + 4z^{16}
\label{e2.4}
\end{array}
\qqq

\nsection{Experimental results}
\label{s3}
In this section we will present the results of numerical calculations
of the coefficients $\dnm(\cK)$ of \eex{1.11},\rx{1.12} for some
simple knots.
L.~Kauffman and S.~Lins\cx{KL} presented conveniently normalized
formulas for the Jones polynomial $\vakk$ as an element of
$\ZZ[\spq]$.
Choose an integer number $N\geq 0$. If
we calculate the polynomial for $1\leq \a \leq N+1$, then
after substituting $\spq = 1 + \xh$ and extracting the coefficients
in front of the first powers $\xh^n$, $0\leq n\leq 2N$ we can
determine all the coefficients $\DmnK$, $0\leq m \leq N$,
$0\leq n\leq 2N$ of \ex{1.4} by solving the system of linear
equations for every power of $\xh$ (each equation in a
particular system corresponds to a specific value of $\a$). Finally
we can recalculate the coefficients $D_{m,n}$ into the coefficients
$\dnm$, $0\leq 2N$, $0\leq m\leq N-{n\over 2}$ of \ex{1.12}. The
results are presented in Tables~\ref{tb.1} and \ref{tb.2} for the
knots $5_2$ and $6_1$ respectively (see \eg \cx{BZ} for the table of
knots). As we see, all the coefficients in the tables are indeed
integer in agreement with the \wc. It is easy to check that the
coefficients in the top line are consistent with the claim\rx{1.14}
of the Melvin-Morton conjecture.

We checked the \sc by using all available coefficients $\dnm$ in
order to calculate `approximate' polynomials $\pappr{n}$ by the
formula
\qq
\pappr{n}(\cK;z) =
\ankz \sum_{0\leq m\leq N-{n\over 2} } \dnm z^{2m} +
\cO(z^{2N-n+1}).
\qqq
The results seem to confirm our conjecture, because the degree of
$\pappr{n}$ appears to be limited:
\qq
\begin{array}{lrcl}
5_2: \qquad &
\nabla_A & = & 1 + 2z^2 \hspace{298.7pt}
\\
& \pappr{1}& = & -6z^2 - 5z^4 + \cO(z^{18})
\\
& \pappr{2}& = &  2 - 7z^2 + 36z^4 + 54z^6 + 23z^8 + \cO(z^{18})
\\
& \pappr{3}& = & 4 - 83z^2 + 140z^4 - 156z^6 - 467z^8 - 358z^{10} -
103z^{12} + \cO(z^{16})

\end{array}
\label{e2.5}
\qqq
\qq
\begin{array}{lrcl}
6_1: \qquad &
\nabla_A & = & 1 - 2z^2 \hspace{298.7pt}
\\
& \pappr{1}& = & 2z^2 - z^4 + \cO(z^{20})
\\
& \pappr{2}& = & -2 + z^2 + 17z^4 - 10z^6 + 3z^8 + \cO(z^{20})
\\
& \pappr{3}& = & -35z^2 +35z^4 + 166z^6 - 113z^8 + 50z^{10} -
11z^{12} + \cO(z^{18})
\end{array}
\label{e2.6}
\qqq

Tables~\ref{tb.3} and \ref{tb.4} contain the lists of the
coefficients $\tdm{2n}$ for the simplest amphicheiral knots $4_1$
(the `figure 8' knot) and $8_3$. All coefficients appear to be
integer in agreement with Corollary~\ref{cor3.1}.
The approximate line polynomials
$\tpappr{2n}$ for the same amphicheiral knots
were calculated by the formula
\qq
\tpappr{2n}(\cK;z) = \nabla_A^{3n+1}(\cK;z) \sum_{m=0}^{N-n}
\tdm{2n} z^{2m} + \cO(z^{2N-2n+1}).
\label{3.9}
\qqq
As we see, they are also of a limited degree:
\qq
\begin{array}{lrcl}
4_1: \qquad &
\nabla_A & = & 1 - z^2 \hspace{298.7pt}
\\
& \tpappr{2}& = & -1 - z^2 + \cO(z^{20})
\\
& \tpappr{4}& = & 4 + 20z^2 + 14z^4 + 2z^6 + \cO(z^{18})
\\
& \tpappr{6}& = & -35 - 430z^2 - 989z^4 - 635z^6 - 140z^8 - 11z^{10}
+ \cO(z^{16})
\end{array}
\label{e2.7}
\qqq
\qq
\begin{array}{lrcl}
8_3: \qquad &
\nabla_A & = & 1 - 4z^2 \hspace{298.7pt}
\\
& \tpappr{2}& = & -4 - 12z^2 + 11z^4 - 4z^6 + \cO(z^{20})
\\
& \tpappr{4}& = & 60 + 1066z^2 + 1482z^4 + 928z^6 + 513z^8
-248z^{10} + 80z^{12} + \cO(z^{18})
\end{array}
\label{e2.8}
\qqq
This confirms the
Conjecture~\ref{c3.2}.

\begin{table}
\qq
\begin{array}{l|*{7} {r} }
m & 0 & 1 & 2 & 3 & 4 & 5 & 6  \\
\hline
d^{(0)}_m &1 &-2   &4     &-8      &16      &-32       &64
\\
d^{(1)}_m &0 &-6   &31    &-114    &360     &-1040     &2832
\\
d^{(2)}_m &2 &-27  &226   &-1286   &5843    &-22974    &81684
\\
d^{(3)}_m &4 &-139 &1750  &-14100  &86613   &-443388   &1991453
\\
d^{(4)}_m &19&-832 &14664 &-158554 &1262646 &-8145921  &45047755
\\
d^{(5)}_m &93&-5720&133890&-1866899&18679183&-148104718&988048870
\end{array}
\nonumber
\qqq
\caption{The coefficients $d^{(n)}_m$ for the knot $5_2$}
\label{tb.1}
\end{table}
\begin{table}
\qq
\begin{array}{l|*{7} {r} }
m & 0 & 1 & 2 & 3 & 4 & 5 & 6 \\
\hline
d^{(0)}_m & 1  & 2    & 4     & 8      & 16      & 32   & 64
\\
d^{(1)}_m & 0  & 2    & 11    & 42     & 136     & 400  &1104
\\
d^{(2)}_m & -2 & -19  & -93   & -340   & -1037   & -2754&-6428
\\
d^{(3)}_m & 0  & -35  & -455  & -3264  & -17389  & -7720&-300255
\\
d^{(4)}_m & 15 & 328  & 2843  & 14830  & 50071   & 74117&-399260
\\
d^{(5)}_m & 13 & 1226 & 24996 & 274355 & 2107672 &12766200&65058967
\end{array}
\nonumber
\qqq
\caption{The coefficients $d^{(n)}_m$ for the knot $6_1$}
\label{tb.2}
\end{table}
%
\begin{table}
\qq
\begin{array}{l|*{7} {r} }
m & 0 & 1 & 2 & 3 & 4 & 5 & 6 \\
\hline
d^{(0)}_m &1 &1    &1     &1       &1       &1         &1
\\
d^{(2)}_m &-1&-5   &-14   &-30     &-55     &-91       &-140
\\
d^{(4)}_m &4 &48   &266   &996     &2926    &7280      &16044
\\
d^{(6)}_m &-35 &-780 &-7214 &-41875 &-180510 &-631436  &-1890680
\\
d^{(8)}_m &543&19434&270472&2251006&13395371&62736271  &245214729
\end{array}
\nonumber
\qqq
\caption{The coefficients $d^{(2n)}_m$ for the knot $4_1$}
\label{tb.3}
\end{table}
%
\begin{table}
\qq
\begin{array}{l|*{7} {r} }
m & 0 & 1 & 2 & 3 & 4 & 5 & 6 \\
\hline
d^{(0)}_m &1 &4    &16    &32      &64      &128       &256
\\
d^{(2)}_m &-4&-76 &-821 &-6868  &-49504  &-323456 &-1970944
\\
d^{(4)}_m &60 &2746 &58210 &840696 &9594881 &93259044 &806300400
%
\end{array}
\nonumber
\qqq
\caption{The coefficients $d^{(2n)}_m$ for the knot $8_3$}
\label{tb.4}
\end{table}

\nsection{Discussion}
\label{s4}

Let us speculate briefly about the possible origin of the \sc and
Conjecture~\ref{c3.2}. We plan to present the expansion\rx{1.8} of
the colored Jones polynomial $\vakk$ of a knot $\cK\in S^3$ as a
contribution of a particular stationary phase point into a certain
finite dimensional integral (such a representation seems to appear
naturally when one tries to derive the expansion\rx{1.8} and \sc from
the universal $R$-matrix). The determinant of the quadratic
form of second derivatives of the exponent in the rapidly oscillating
exponential is equal to the Alexander polynomial of $\cK$. Thus,
in accordance with \ex{1.14}, the stationary phase point contribution
is inversely proportional to the Alexander polynomial in the
leading approximation in $K^{-1}$.

The subleading terms in the $K^{-1}$ expansion of $\vakk$ can be
calculated by Feynman rules. In other words, the lines $\vnkz$ of the
expansion\rx{1.11} will be related to closed $(n+1)$-loop graphs.
The edges of the graphs represent the inverse matrix of the second
derivative quadratic form. The valence of the vertices of the graphs
matches the order of the higher order terms in the Taylor expansion
of the rapidly oscillating exponent around the stationary phase
point.

The matrix elements of the inverse quadratic form are inversely
proportional to the determinant of that form which is equal to the
Alexander polynomial. As a result, the highest order of the Alexander
polynomial in the denominator of $\vnkz$ will be equal to the maximum
number of edges in a closed $(n+1)$-loop diagram (plus 1 coming from
the leading approximation). For a fixed $n$, this number is
determined by vertices with the smallest valence $\vmin$:
\qq
\#{\rm edges} \leq {\vmin \over \vmin - 2}.
\label{4.1}
\qqq
The exponent of the rapidly oscillating exponential turns out
to be even. Hence, for a general knot, $\vmin=4$, so that the \rhs of
\ex{4.1} is equal to $2n$. This leads to the power $2n+1$ in the
denominator of \ex{1.17}.

One may speculate that for the amphicheiral knots the 4-valent
vertices are absent. Then $\vmin=6$, the \rhs of \ex{4.1} is equal
to ${3\over 2}$, and we reproduce the power in the denominator of
\ex{3.8}.

\hyphenation{Thur-ston}

\section*{Acknowledgements}
I am thankful to D.~Bar-Natan, J.~Bernstein, I.~Cherednik,
P.~Deligne, G.~Felder, D.~Freed, S.~Garoufalidis, J.~Handfield,
L.~Kauffman, M.~Kontsevich, P.~Melvin, J.~Roberts, D.~Thurston,
V.~Turaev, A.~Vaintrob and E.Witten for many useful discussions.

This work was supported by the National Science Foundation
under Grant DMS 9304580. Part of this work was
done during my visit to Argonne National Laboratory. I want to thank
Professor C.~Zachos for his hospitality and advice.

\pagebreak
\section*{Figure Captions}

\begin{description}
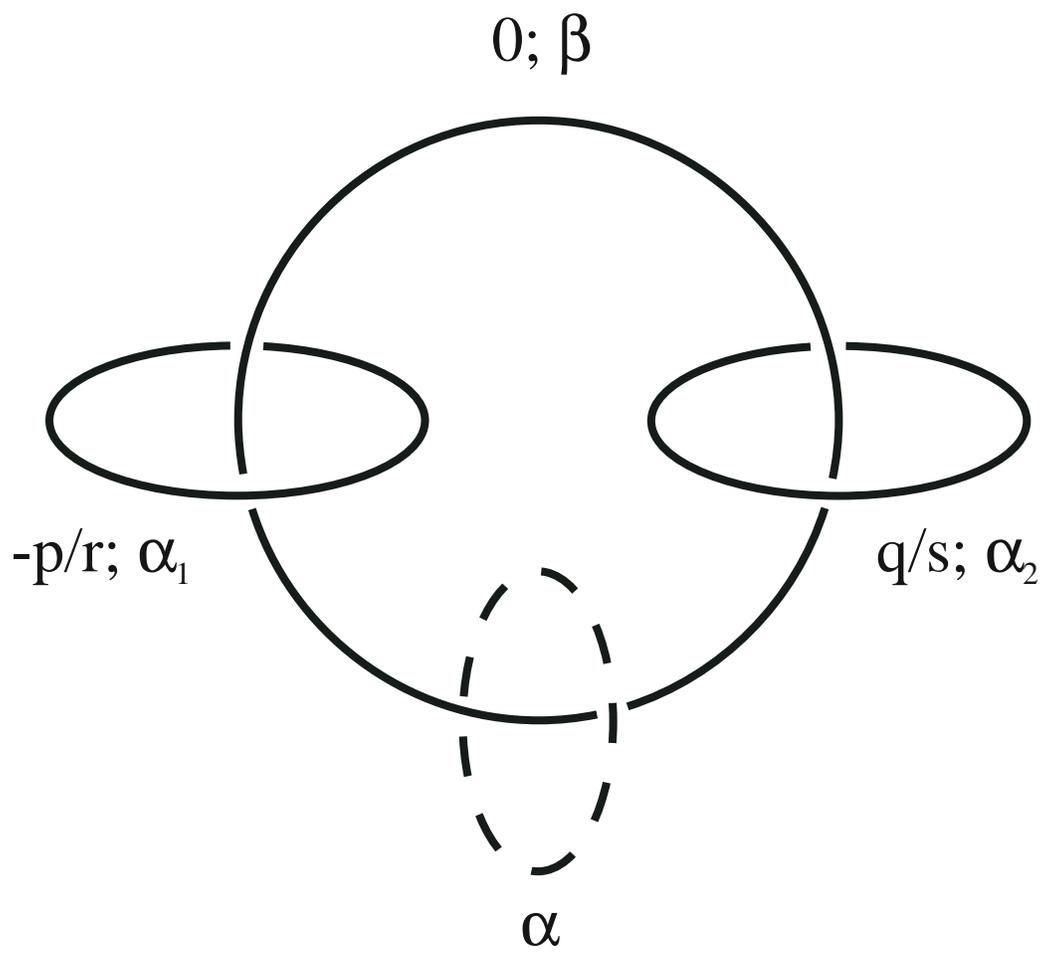
\item[Fig.~1]
A surgery link for producing a torus knot from an unknot.

\end{description}
\pagebreak

\end{document}

\end{document}

The equation\rx{2.19} would have proved the \sc if not for the
fractions in the expansion of
$(1+\xh)^{ {1\over 4}\pqf }$ and
$\log(1+\xh)\over 4\xp\xq$ in powers of $\xh$. Thus we have proved so
far only that

Our strategy is to show that
$\pnpq \in \zzf{\xp}$ and similarly $\pnpq \in \zzf{\xq}$. Since
$\xp$ and $\xq$ are coprime, $\zzf{\xp} \bigcap \zzf{\xq}=\zzs{z^2}$
and
this will prove the \sc.

\end{document}